# STRONG GAUSSIAN APPROXIMATIONS OF PRODUCT-LIMIT AND QUANTILE PROCESSES FOR STRONG MIXING AND CENSORED DATA

V. FAKOOR      AND      N. NAKHAEE RAD

*Ferdowsi University of Mashhad, and Azad University of Mashhad*

ABSTRACT. In this paper, we consider the product-limit quantile estimator of an unknown quantile function under a censored dependent model. This is a parallel problem to the estimation of the unknown distribution function by the product-limit estimator under the same model. Simultaneous strong Gaussian approximations of the product-limit process and product-limit quantile process are constructed with rate $O((\log n)^{-\lambda})$ for some $\lambda > 0,$. The strong Gaussian approximation of the product-limit process is then applied to derive the laws of the iterated logarithm for product-limit process.

## 1. Introduction and Preliminaries

In medical follow-up or in engineering life testing studies, the life time variable may not be observable. Let $X_1, \ldots, X_n$ be a sequence of life times, having a common unknown continuous marginal distribution function $F$, with a density function $f = F'$ and hazard rate $\lambda = f/(1 - F)$. The random variables are not assumed to be mutually independent (see Assumption (1) for the kind of dependence stipulated). Let the random variable $X_i$ be censored on the right by the random variable $Y_i$, so that one observes only

$$Z_i = X_i \wedge Y_i \qquad \text{and} \qquad \delta_i = I(X_i \leq Y_i),$$

where $\wedge$ denotes minimum and $I(.)$ is the indicator of the event specified in parentheses. In this random censorship model, we assume that the censoring random variables $Y_1, \ldots, Y_n$ are not mutually independent (see Assumption (2) for the kind of dependence stipulated), having a common unknown continuous d.f. $G$, and that they are independent of the $X_i$'s. Since censored data traditionally occur in lifetime analysis, we assume that $X_i$ and $Y_i$ are nonnegative. The actually observed $Z_i$'s have a distribution function $H$ satisfying

$$\overline{H}(t) = 1 - H(t) = (1 - F(t))(1 - G(t)).$$

Denote by

$$F_*(t) = P(Z \leq t, \delta = 1),$$

the sub-distribution function for the uncensored observations. Define

$$N_n(t) = \sum_{i=1}^{n} I(Z_i \leq t, \delta = 1) = \sum_{i=1}^{n} I(X_i \leq t \wedge Y_i),$$

the number of uncensored observations less than or equal to $t$, and

$$Y_n(t) = \sum_{i=1}^{n} I(Z_i \geq t),$$

the number of censored or uncensored observations greater than or equal to $t$ and also the empirical distribution functions of $\bar{H}(t)$ and $F_*(t)$ are respectively defined as

$$\overline{Y}_n(t) = n^{-1}Y_n(t) \quad , \quad \overline{N}_n(t) = n^{-1}N_n(t).$$

Then the Kaplan-Meier estimator for $1 - F(t)$, based on the censored data is

$$(1.1) \qquad 1 - \widehat{F}_n(t) = \prod_{s \leq t}(1 - \frac{dN_n(s)}{Y_n(s)}), \qquad t < Z_{(n)},$$

where $Z_{(i)}$ are the order statistics of $Z_i$ and $dN_n(t) = N_n(t) - N_n(t-)$.

As is known (see, e.g., Gill, 1980), for a d.f. $F$ on $[0, \infty)$, the cumulative hazard function $\Lambda$ is defined by

$$\Lambda(t) = \int_0^t \frac{dF(s)}{1 - F(s^-)},$$

and $\Lambda(t) = -log(1 - F(t))$ for the case that $F$ is continuous. The empirical cumulative hazard function $\hat{\Lambda}_n(t)$ is given by

$$\widehat{\Lambda}_n(t) = \int_0^t \frac{dN_n(s)}{Y_n(s)},$$

which is called the Nelson-Aalen estimator of $\Lambda(t)$ in the literature.

For a censored model with $\{X_i; i \geq 1\}$ and $\{Y_i; i \geq 1\}$ independent and identically distributed (i.i.d) sequences and mutually independent, Burke et al. (1981, 1988), established strong Gaussian approximation of the product-limit(PL) process $\sqrt{n}(\widehat{F}_n(t) - F(t))$ by a two parameter Gaussian process at the almost sure rate of $O(n^{-1/2}(\log n)^2)$. In left truncation and right censorship(LTRC) model, Zhou and Yip (1999) initiated and Tse (2003, 2005) established strong Gaussian approximation of the PL-process by a two-parameter Gaussian process at the almost sure rate of $O((\log n)^{3/2}n^{-1/8})$, a rate that reflects the two-dimensional nature of the LTRC model.

The quantile function $Q$ and its sample estimator $Q_n$ are defined by

$$Q(p) := \inf\{t : F(t) \geq p\}, \qquad Q_n(p) := \inf\{t : \widehat{F}_n(t) \geq p\}$$

for $0 < p < 1$. The role of the quantile function in statistical data modeling was emphasized by Parzen (1979). In econometrics, Gastwirth (1971) used the quantile function to give a succinct definition of the Lorenz curve, which measure inequality in distribution of resources and in size distribution.

In the independent framework with no censoring, the properties of estimator $Q_n$ (where $\widehat{F}_n$ is replaced by the empirical d.f. $F_n$) have been extensively studied (see e.g., Csörgő, 1983; Shorack and Wellner, 1986). Under $\phi$-mixing condition (for the definition see

Doukhan, 1996), the Bahadur representation was obtained by Sen (1972) and the extension to the $\alpha$-mixing case was obtained by Yoshihara (1995). Under $\alpha$-mixing condition(see definition below), the strong approximation of the normed PL-quantile process $\rho_n(p) := \sqrt{n} f(Q(p))[Q(p) - Q_n(p)]$ by a two parameter Guassian process at the rate $O((\log n)^{-\lambda})$ for some $\lambda > 0$, was obtained by Fotopoulos et al.(1994) and was later improved by Yu (1996).

For a censored model with $X_i$ and $Y_i$'s, independent and identically distributed sequences and mutually independent, Padgett and Thombs (1989) stated the strong consistency and asymptotic normality for a smooth estimator of $Q(p)$. Sander (1975) obtained some asymptotic properties, and Csörgő (1983) and Cheng (1984) discussed strong approximation results with some applications for $Q_n(p)$. In left truncation and right censorship model, Tse (2005), obtained strong Gaussian approximations of the PL-quantile process by a two parameter Kifere type process at the rate $O((\log n)^{3/2} n^{-1/8})$. Ould-Saïd and Sadki (2005) established the strong consistency and a Badadur-type representation of K-M quantile function $Q_n(.)$ under a strong mixing hypothesis.

The main aim of this paper is to derive strong Gaussian approximations of the PL-process and PL-quantile process, for the case in which the underling lifetime are assumed to be $\alpha$-mixing whose definition is given below. As a result, we obtain the Law of the iterated logarithm for PL-process.

For easy reference, let us recall the following definition.

**Definition 1.** *Let $\{X_i, i \geq 1\}$ denote a sequence of random variables. Given a positive integer $m$, set*

(1.2) $$\alpha(m) = \sup_{k \geq 1} \{|P(A \cap B) - P(A)P(B)| \; ; \; A \in \mathcal{F}_1^k, B \in \mathcal{F}_{k+m}^\infty\},$$

*where $\mathcal{F}_i^k$ denote the $\sigma$-field of events generated by $\{X_j; i \leq j \leq k\}$. The sequence is said to be $\alpha$-mixing (strongly mixing) if the mixing coefficient $\alpha(m) \to 0$ as $m \to \infty$.*

Among various mixing conditions used in the literature, $\alpha$-mixing, is reasonably weak and has many practical applications. There exists many processes and time series fulfilling the strong mixing condition. In particular, the stationary autoregressive-moving average (ARMA) processes, which are widely applied in time series analysis, are $\alpha$-mixing with exponential mixing coefficient, i.e., $\alpha(n) = e^{-\nu n}$ for some $\nu > 0$. The threshold models, the EXPAR models(see Ozaki, 1979), the simple ARCH models(see Engle, 1982; Masry and Tjostheim, 1995, 1997) and their extensions(see Diebolt and Guégan, 1993) and the bilinear Markovian models are geometrically strongly mixing under some general ergodicity conditions. Auestad and Tjostheim (1990) provided excellent discussions on the role of $\alpha$-mixing for model identification in nonlinear time series analysis.

Now, for the sake of simplicity, the assumptions used in this paper are as follows.

**Assumptions.**

(1) Suppose that $\{X_i, i \geq 1\}$ is a sequence of stationary $\alpha$-mixing random variables with continuous distribution function $F$ and mixing coefficient $\alpha_1(n)$.

(2) Suppose that $\{Y_i, i \geq 1\}$ is a sequence of stationary $\alpha$-mixing random variables with continuous distribution function $G$ and and mixing coefficient $\alpha_2(n)$. Moreover, we assume the censoring times are independent of $\{X_i, i \geq 1\}$.

(3) $\alpha(n) = O(e^{-\log^{1+\zeta} n})$ for some $\zeta > 0$, with $\alpha(n) = \max(\alpha_1(n), \alpha_2(n))$ (see Remark 2.1. in Ould-Saïd and Sadki (2005)).

Cai(Lemma 1, 1998) showed that the $Z_i$'s are $\alpha$-mixing random variables(with appropriate coefficient $\alpha$).

The layout of the paper is as follows. Section 2, contains main results. The proofs of the main results are relegated to Section 3.

## 2. Main Results

In the first theorem, we construct a two parameter mean zero Gaussian process that strongly uniformly approximate the empirical processes $Z_{n1}(t) = \sqrt{n}(\widehat{\Lambda}_n(t) - \Lambda(t))$ and $Z_{n2}(t) = \sqrt{n}(\widehat{F}_n(t) - F(t))$.

**Theorem 1.** *Suppose that Assumptions (1)-(3) is satisfied. On a rich probability space, there exists a two parameter mean zero Gaussian process $\{B(u,v)\ u, v \geq 0\}$ such that,*

$$(2.1) \qquad \sup_{t \geq 0} |Z_{n1}(t) - B(t,n)| = O((log n)^{-\lambda}) \quad a.s.,$$

$$(2.2) \qquad \sup_{t \geq 0} |Z_{n2}(t) - (1 - F(t))B(t,n)| = O((log n)^{-\lambda}) \quad a.s.,$$

*for $\lambda > 0$.*

**Remark 1.** *In the $\alpha$-mixing case, we can not achieve the same rate as in the iid case i.e. $O(n^{-1/2}(\log n)^2)$ (see Burke et al.(1988), Theorem 1). The main reason is that our approach utilizes the strong approximation introduced by Dhompongsa(1984) as a kiefer process with a negligible reminder term of order $O(n^{-1/2}(\log n)^{-\lambda})$. This is not as sharp as in iid case.*

**Corollary 1.** *Under assumptions $(1) - (3)$, we have,*

$$(2.3) \qquad \sup_{t \in R} |\widehat{\Lambda}_n(t) - \Lambda(t)| = O\left(\frac{\log \log n}{n}\right)^{1/2} \quad a.s.,$$

$$(2.4) \qquad \sup_{t \in R} |\widehat{F}_n(t) - F(t)| = O\left(\frac{\log \log n}{n}\right)^{1/2} \quad a.s.$$

In the next theorem, we construct a two parameter mean zero Gaussian process that strongly uniformly approximate the empirical process $\rho_n(p)$.

**Theorem 2.** *Let $0 < p_0 \leq p_1 < 1$. Under assumptions (1)-(3), assume that $F$ is Lipschtiz continuous and that $F$ is twice continuously differentiable on $[Q(p_0) - \delta, Q(p_1) + \delta]$ for some* $\delta > 0$ such that $f$ is bounded away from zero. Then there exists a two parameter mean zero Gaussian process $B(t, u)$ for $t, u \geq 0$ such that,

$$(2.5) \qquad \sup_{p_0 \leq p \leq p_1} |\rho_n(p) - (1-p)B(Q(p), n)| = O((logn)^{-\lambda}) \quad a.s.,$$

with $\lambda > 0$.

## 3. Proofs

In order to prove Theorem 1, we need the following lemma.

**Lemma 1.** *(Theorem 3 in Dhompongsa 1984). Under assumptions (1) and (3), there exists a Kiefer process $\{K(s, t), \ s \in \mathbb{R}, t \geq 0\}$ with covariance function*

$$E[K(s,t)K(s',t')] = \Gamma(s,s')\min(t,t')$$

and $\Gamma(s, s')$ is defined by

$$\Gamma(s,s') = Cov(g_1(s), g_1(s')) + \sum_{k=2}^{\infty}[Cov(g_1(s), g_k(s')) + Cov(g_1(s'), g_k(s))],$$

*where $g_k(s) = I(Z_k \leq s) - H(s)$, such that, for some $\lambda > 0$ depending only on $\nu$, given in assumption (3),*

$$\sup_{t \in \mathbb{R}} |\overline{Y}_n(t) - \overline{H}(t) - K(t, n)/n| = O(b_n), \quad a.s.$$

*where*

$$b_n = n^{-1/2}(\log n)^{-\lambda}.$$

**Proof of Theorem 1.** We start with the usual decomposition of $Z_{n1}(t)$.

$$
\begin{aligned}
Z_{n1}(t) &= \sqrt{n}[\hat{\Lambda}_n(t) - \Lambda(t)] = \sqrt{n}\left[\int_0^t \frac{d\overline{N}_n(x)}{\overline{Y}_n(x)} - \int_0^t \frac{dF_*(x)}{\overline{H}(x)}\right] \\
&= \int_0^t \frac{\sqrt{n}[\overline{H}(x) - \overline{Y}_n(x)]}{(\overline{H}(x))^2}dF_*(x) + R_{n1}(t)
\end{aligned}
$$

where

$$
\begin{aligned}
n^{-1/2}R_{n1}(t) &= \int_0^t \frac{(\overline{Y}_n(x) - \overline{H}(x))^2}{\overline{Y}_n(x)(\overline{H}(x))^2}dF_*(x) + \int_0^t \frac{d[\overline{N}_n(x) - F_*(x)]}{\overline{H}(x)} \\
&+ \int_0^t \left(\frac{1}{\overline{Y}_n(x)} - \frac{1}{\overline{H}(x)}\right)d[\overline{N}_n(x) - F_*(x)] \\
&= I_1 + I_2 + I_3.
\end{aligned}
$$

Define, for $t \geq 0$ the sequence of Gaussian processes

$$(3.1) \qquad B(t, n) = \int_0^t \frac{K(x, n)/\sqrt{n}}{(\overline{H}(x))^2}dF_*(x),$$

where $K(s,t)$ is the Kiefer process in Lemma 1. Clearly, $E(B(t,n)) = 0$,

$$Cov[B(s,m), B(t,n)] = \sqrt{\frac{m}{n}} \int_0^t \int_0^s \frac{\Gamma(x,y)}{\overline{H}(x)^2 \overline{H}(y)^2} dF_*(x) dF_*(y).$$

where $\Gamma(s,t)$ is defined in Lemma 1. Let

$$\beta(t,n) = -\sqrt{n}(\overline{Y}_n(t) - \overline{H}(t)) - K(t,n)/\sqrt{n}.$$

Theorem 1 is about the order

$$(3.2) \qquad \sup_{t \geq 0} |\widehat{Z}_{n1}(t) - B(t,n)| = \sup_{t \geq 0} |R_{n1}(t) + R_{n2}(t)|,$$

where

$$R_{n2}(t) = \int_0^t \frac{\beta(x,n)}{(\overline{H}(x))^2} dF_*(x).$$

To deal with $R_{n1}(t)$, we treat $\overline{Y}_n(x)$ as an empirical d.f. associated to $Z_i's$ and from Theorem 1 in Cai (1998), we have

$$(3.3) \qquad I_1 = O(a_n^2) \quad a.s.,$$

where

$$a_n = \left(\frac{\log\log n}{n}\right)^{1/2}.$$

To estimates of $I_2$, divide the interval $[0,t]$ into subintervals $[x_i, x_{i+1}], i = 1, \ldots, k_n$ where $k_n = O(a_n^{-1})$, and $0 = x_1 < x_2 <, \cdots < x_{k_n+1} = t$ are such that $x_{i+1} - x_i = O(b_n)$. Then

$$|I_2| = \left|\int_0^t \frac{d[\overline{N}_n(x) - F_*(x)]}{\overline{H}(x)}\right| \leq \sum_{i=0}^{k_n+1} \left|\int_{x_i}^{x_{i+1}} \frac{d[\overline{N}_n(x) - F_*(x)]}{\overline{H}(x)}\right|.$$

The integral on the right hand side of the latter inequality is bounded above by $(\log n)^{-\lambda-\beta}/\sqrt{n}$, almost surely(The proof of this can be done using similar arguments to $A$ in Lemma 3.4, in Ould-Saïd et al. 2005). Therefore

$$(3.4) \qquad I_2 = O((\log n)^{-\lambda-\beta}) \quad a.s.$$

The estimate of $I_3$ is similar to the estimate of $I_3$ in the proof of Theorem 2 in Cai(1998). Hence,

$$(3.5) \qquad I_3 = O(b_n) \quad a.s.$$

Therefore, by combining (3.3)-(3.5), we have

$$(3.6) \qquad \sup_{t \geq 0} |R_{n1}(t)| = O((\log n)^{-\lambda-\beta}) \quad a.s.$$

Next, by applying Lemma 1, we have

$$(3.7) \qquad \sup_{t \geq 0} |R_{n2}(t)| = O((\log n)^{-\lambda}) \quad a.s.$$

Combining (3.2), (3.6) and (3.7) we obtain (2.1). It can be shown that

$$(3.8) \qquad \widehat{F}_n(t) - F(t) = (1 - F(t))[\hat{\Lambda}_n(t) - \Lambda(t)] + O(\frac{\log\log n}{n}) \quad a.s.$$

Therefore (2.2) is proved via (3.8). $\square$

**Proof of Corollary 1.** By the law of the iterated logarithm for Kiefer processes (see, e.g., Corollary 1.15.1, page 81, in Csörgő and Révész, 1981), and (3.1) we have,

$$(3.9) \qquad \sup_{t \in R} |B(t,n)| \leq C \sup_{t \in R} |K(t,n)|/\sqrt{n} = O\left(\frac{\log \log n}{n}\right)^{1/2} \quad a.s.$$

where $C$ is a positive constant. From (2.1), (2.2) and (3.9) we obtain the results. $\square$

The proof of Theorem 2, is mainly based on the following Lemmas of Ould-Saïd et. al. (2005). Lemma 2 shows that $\widehat{F}_n$ composed with $Q_n$ is an approximate identity up to order $O(b_n)$. Lemmas 3 and 4 give global and local bounds for the deviation between $Q_n$ and $Q$.

**Lemma 2.** *Let $0 < p_0 \leq p_1 < 1$. Under assumptions (1)-(3), assuming that $F$ is continuous, then*

$$\sup_{p_0 \leq p \leq p_1} |\widehat{F}_n(Q_n(p)) - p| = O(b_n) \quad a.s.$$

*where $b_n = (\log n)^{-\lambda}/\sqrt{n}$ with $\lambda > 0$.*

**Lemma 3.** *Under assumptions (1)-(3), assuming that $F' = f$ is continuous and strictly positive on $[Q(p_0) - \delta, Q(p_1) + \delta]$ for some $\delta > 0$. Then,*

$$\sup_{p_0 \leq p \leq p_1} \sqrt{n}|Q_n(p) - Q(p)| = O(\sqrt{\log \log n}) \quad a.s.$$

**Lemma 4.** *Let $\lambda_n = const.b_n$. Under (1)-(3) and assuming that $F$ is Lipschitz continuous on $[0, \tau]$. Then,*

$$\sup_{|s-t| \leq \lambda_n} \sqrt{n}|Z_{n2}(s) - Z_{n2}(t)| = O((\log n)^{-\beta-\lambda}) \quad a.s.,$$

*where $\beta > 0$ and $\lambda > 0$.*

**Proof of Theorem 2.** We continue to use the notation $\lambda_n$ as in Lemma 4. Let $s = Q_n(p)$ and $t = Q(p)$, $p_0 \leq p \leq p_1$, Lemma 3 yields $\sqrt{n}|s - t| = O(\sqrt{\log \log n})$. Applying Lemma 4 gives,

$$(3.10) \qquad \widehat{F}_n(Q_n(p)) - \widehat{F}_n(Q(p)) = F(Q_n(p)) - F(Q(p)) + O((\log n)^{-\beta-\lambda}/\sqrt{n}) \quad a.s.$$

By Lemma 2, $\widehat{F}_n(Q_n(p))$ can be replaced by $p$ up to $O(b_n)$. For the right hand side, a Taylor expansion of the first term about $Q(p)$ up to second order term gives,

$$f(Q(p))[Q(p) - Q(p)] + O([Q_n(p) - Q(p)]^2) + O((\log n)^{-\beta-\lambda}/\sqrt{n}) \quad a.s.,$$

for $p_0 \leq p \leq p_1$. Invoking Lemma 3 and rearranging terms in (3.10), we have,

$$\sqrt{n}f(Q(p))[Q_n(p) - Q(p)] = \sqrt{n}[p - \widehat{F}_n(Q(p))] + O((\log n)^{-\beta-\lambda}) \quad a.s.,$$

for $p_0 \leq p \leq p_1$. Since $F$ is continuous, $F(Q(p)) = p$. Recalling the definitions of the PL process $Z_n$ and PL-quantile process $\rho_n$, we have,

$$(3.11) \qquad \rho_n(p) = \widehat{Z}_n(Q(p)) + O((\log n)^{-\beta-\lambda}/\sqrt{n}) \quad a.s.$$

for $p_0 \leq p \leq p_1$. By using Theorem 1 and (3.11), theorem is proved. $\square$

## Acknowledgements

The first author wish to acknowledge partial support from Statistics Center of Excellence of Ferdowsi University of Mashhad

Department of Statistics,
School of Mathematical Sciences,
Ferdowsi University of Mashhad, Iran.
P. O. Box: 1159-91775.

Department of Statistics,
Azad University of Mashhad, Iran.

*E-mail address*: `fakoor@math.um.ac.ir`
*E-mail address*: `najme.rad@gmail.com`